\documentclass[runningheads]{llncs}
\usepackage[T1]{fontenc}
\usepackage[utf8]{inputenc}
\usepackage{amsmath,amssymb}
\usepackage{isabelle,isabellesym}
\usepackage{graphicx}

\usepackage{hyperref,color}

\begin{document}
\title{Large-Scale Formal Proof for the Working Mathematician --- Lessons learnt from the ALEXANDRIA Project}
\titlerunning{Large-Scale Formal Proof --- Lessons from ALEXANDRIA}
\author{Lawrence C Paulson\orcidID{0000-0003-0288-4279}}
\authorrunning{L C Paulson}
\institute{Computer Laboratory, University of Cambridge, UK\\
\email{lp15@cam.ac.uk} \quad
\url{https://www.cl.cam.ac.uk/~lp15/}}
%
%
%
\maketitle              
\begin{abstract}
ALEXANDRIA is an ERC-funded project that started in 2017, with the aim of bringing formal verification to mathematics. The past six years have seen great strides in the formalisation of mathematics and also in some relevant technologies, above all machine learning. Six years of intensive formalisation activity seem to show that even the most advanced results, drawing on multiple fields of mathematics, can be formalised using the tools available today.

\end{abstract}

\section{Introduction}

In the summer of 2017, the Newton Institute at Cambridge held a programme entitled \emph{Big Proof} (BPR) ``directed at the challenges of bringing proof technology into mainstream mathematical practice''. It was held in recognition of the formalisations that had already been done (which were indeed big). The programme webpage%
\footnote{\url{https://www.newton.ac.uk/event/bpr/}}
specifically lists the proofs of the Kepler conjecture~\cite{hales-formal-Kepler}, the odd order theorem~\cite{gonthier-oot} and the four colour theorem~\cite{gonthier-4ct}. That summer also saw the start of my ERC project, ALEXANDRIA. \emph{Big Proof} represented an acknowledgement that the formalisation of mathematics could no longer be ignored, but also an assertion that big problems remain to be solved.These included ``novel pragmatic foundations'' and large-scale ``formal mathematical libraries'' and ``inference engines'', and also the ``curation'' of formalised mathematical knowledge.

ALEXANDRIA was conceived in part to try to identify those big problems. By hiring professional mathematicians and asking them to formalise advanced mathematics, we would get a direct idea of the obstacles they faced. We would also try to refine our tools, extend our libraries and investigate other technologies. We would have only five years (extended to six due to COVID-19).

The need for formalisation had been stressed by Vladimir Voevodsky, a Fields medallist, who pointedly asked ``And who would ensure that I did not forget something and did not make a mistake, if even the mistakes in much more simple arguments take years to uncover?'' \cite{voevodsky-origins}. He advocated a new sort of formalism, homotopy type theory, which was the subject of much excitement. However, the most impressive formalisations by that time had been done in Coq (four colour theorem, odd order theorem), HOL Light (Kepler conjecture and much else) or Isabelle/HOL (part of the Kepler proof, and more). Lean, a newcomer, was attracting a user community. Perhaps our project would shed light on the respective values of the available formalisms: calculus of constructions (Coq, Lean), higher-order logic or homotopy type theory. Voevodsky would never find out, due to his untimely death in September 2017.

Since that date, research into the formalisation of mathematics has plunged ahead. Kevin Buzzard, a number theorist at Imperial College London, followed some of the \emph{Big Proof} talks online. This resulted in his adoption of Lean for his Xena Project, with the aim of attracting students to formalisation.%
\footnote{\url{https://www.ma.imperial.ac.uk/~buzzard/xena/}} 
Xena has had a huge impact, but here I'd like to focus on the work done within ALEXANDRIA\@.

\section{A Brief Prehistory of the Formalisation of Mathematics}

Mathematics is a work of the imagination, and the struggle between intuition and rigour has gone on since classical times. Euclid's great contribution to Greek geometry was the unification of many separate schools through his system of axioms and postulates. Newton and Leibniz revolutionised mathematics, but the introduction of infinitesimals was problematical. During the 19th centuries, the ``arithmetisation of analysis'' carried out by Cauchy and Weierstrass replaced infinitesimals by rigorous $\epsilon$--$\delta$ arguments. (We would not get a consistent theory of infinitesimals until the 1960s, under the banner of non-standard analysis.) Dedekind and Cantor promulgated a radical new understanding of sets and functions that turned out to be inconsistent until Zermelo came up with his axioms. It is notable that Zermelo set theory (which includes the axiom of choice but lacks Fraenkel's replacement axiom) is approximately equal in logical strength to higher-order logic.

Only axiomatic mathematics can be formalised. The first attempt was by Frege, whose work (contrary to common belief) was not significantly impacted by Russell's paradox. Russell and Whitehead in their \emph{Principia Mathematica}~\cite{principia} wrote out the proofs of thousands of mathematical propositions in a detailed axiomatic form. The work of Bourbaki can also be seen as a kind of formalised mathematics. The philosopher Hao Wang wrote on the topic and also coded the first automatic theorem prover~\cite{wang-toward-mechanical} for first-order logic, based on what we would now recognise as a tableau calculus.

This takes us to NG de Bruijn, who in 1968 created AUTOMATH~\cite{debruijn68}, and to his student's formalisation~\cite{jutting77} of Landau's \emph{Foundations of Analysis} in 1977. This takes us to the birth of Mizar~\cite{Grabowski2015}, in which a truly impressive amount of mathematics was formalised in a remarkably readable notation. More recent history --- analysis in HOL Light, the four colour theorem in Coq, etc --- is presumably familiar to readers. But it is appropriate to close this section with a prescient remark by de Bruijn back in 1968:

\begin{quote}
  As to the question what part of mathematics can be written in AUTOMATH, it should first be remarked that we do not possess a workable definition of the word "mathematics". 
Quite often a mathematician jumps from his mathematical language into a kind of metalanguage, obtains results there, and uses these results in his original context. It seems to be very hard to create a single language in which such things can be done without any restriction.\cite[p.\ts3]{AUT001}
\end{quote}

And so we have two great scientific questions: 
\begin{itemize}
  \item \textbf{What sort of mathematics can be formalised?}
  \item \textbf{What sort of proofs can be formalised?}
\end{itemize}

We would investigate these questions --- mostly in the context of Isabelle/HOL --- by formalising as much mathematics as we could, covering as many different topics as possible. I expected to run into obstacles here and there, which would have to be recorded if they could not be overcome.

\section{ALEXANDRIA: Warmup Formalisation Exercises}

The ERC proposal called for hiring research mathematicians, who would bring their knowledge of mathematics as it was practised, along with their \emph{inexperience} of Isabelle/HOL\@. Their role would be to formalise increasingly advanced mathematical material with the twin objectives of developing formalisation methodologies and identifying deficiencies that might be remedied by extending Isabelle/HOL somehow. The project started in September 2017. We hired Anthony Bordg and Angeliki Koutsoukou-Argyraki. A third postdoc was required to undertake any necessary Isabelle engineering, and Wenda Li was hired.

One of the tasks for the first year was simply to reorganise and consolidate the Isabelle/HOL analysis library, which had mostly been translated from HOL Light. But we were also supposed to conduct pilot studies. The team set to work enthusiastically, and already in the first year they created a number of impressive developments:
\begin{itemize}
  \item \emph{Irrational rapidly convergent series}, formalising a 2002 proof by J. Hančl~\cite{hancl}
  \item \emph{Projective geometry}, including Hessenberg's theorem and Desargues's theorem
  \item The theory of \emph{quantum computing} (which identified a significant error in one of the main early papers)
  \item \emph{Quaternions}, \emph{octonions} and several other small exercises
  \item Effectively counting \emph{real and complex roots of polynomials}, and the Budan-Fourier theorem \cite{li-counting-polynomial,li-evaluating-winding}
  \item The first formal proof that \emph{every field contains an algebraically closed extension}~\cite{vilhena-algebraically-closed}
\end{itemize}
Koutsoukou-Argyraki wrote up her reactions to Isabelle/HOL from the perspective of a mathematician in her paper ``Formalising Mathematics --- in Praxis''~\cite{koutsoukou-argyraki-formalising-mathematics}.

\section{Advanced Formalisations}

As noted above, Kevin Buzzard had taken an interest in formalisation through participation in \emph{Big Proof}, and by 2019 had marshalled large numbers of enthusiastic students to formalise mathematics using Lean. He had also made trenchant criticisms of even the most impressive prior achievements: that most of it concerned simple objects such as finite groups, or was just 19th-century mathematics. Nobody seemed to be working with sophisticated objects. He expressed astonishment that Grothendieck schemes --- fundamental objects in algebraic geometry and number theory --- had not been formalised in any tool. His criticisms helped focus our attention on the need to tackle difficult, recent and deep mathematics. Team members proposed their own tasks, but we also contributed to one another's tasks, sometimes with the help of interns or students. We completed three notable projects during this middle period:

\begin{itemize}
  \item \emph{Irrationality and transcendence criteria for infinite series}~\cite{koutsoukou-irrationality}, extending the Hančl work mentioned above with material from two more papers: Erdős--Straus~\cite{erdos-irrationality} and Hančl--Rucki~\cite{hancl2005}.
  \item \emph{Ordinal partition theory}~\cite{dzamonja-formalising}: infinite forms of Ramsey's theorem, but for order types rather than cardinals. We formalised relatively  papers by Erd\H{o}s--Milner~\cite{erdos-theorem-partition} and Larson~\cite{larson-short-proof}, and as a preliminary, the Nash-Williams partition theorem~\cite{todorcevic-ramsey}. These were deep results in the context of Zermelo--Fraenkel set theory, involving highly intricate inductive constructions. One of the papers contained so many errors as to necessitate publishing a second paper~\cite{erdos-theorem-partition-corr} with a substantially different proof. This material was difficult even for Erd\H{o}s!
  \item \emph{Grothendieck Schemes}~\cite{bordg-simple}. Buzzard had formalised schemes in Lean~\cite{buzzard-schemes-lean} (three times), and even claimed that Isabelle was not up to the job due to its simple type system. We took the challenge and it was straightforward, following a new approach based on locales to manage the deep hierarchies of definitions. 
\end{itemize}

 We were aiming for a special issue devoted to formalisation in the journal \emph{Experimental Mathematics}, and were delighted to see these projects take up three of the six papers ultimately accepted. 
 
 \section{Seriously Deep Formalisation Projects}
 
 Inspired by the success of the previous projects --- conducted under the difficult circumstances of COVID-19 lockdown --- team members continued to propose theorems to formalise, and we continued to collaborate in small groups. By now we had the confidence to take on almost anything. There are too many projects to describe in full, so let's look at some of the highlights.
 
\subsection{Szemerédi's regularity lemma and Roth's theorem on arithmetic progressions}
 
\emph{Szemerédi's regularity lemma} is a fundamental result in extremal graph theory. It concerns a property called the \emph{edge density} of two given sets of vertices $X$, $Y\subseteq V(G)$, and a further property of $(X,Y)$ being an $\epsilon$-regular pair for any given $\epsilon>0$. The lemma itself states that for a given $\epsilon>0$ there exists some $M$ such that every graph has an $\epsilon$-regular partition of its vertex set into at most $M$ parts. 
Intuitively, $(X,Y)$ is an $\epsilon$-regular pair if the density of edges between various subsets $A\subseteq X$ and $B\subseteq Y$ is more or less the same for all possible $A$ and~$B$; an $\epsilon$-regular partition enjoys that property for all but an insignificant number of pairs $(X,Y)$ of vertex sets taken from the partition. Intuitively then, the vertices of any graph can be partitioned into most $M$ parts such that the edges between the various parts are uniform in this sense.
 
 We used Szemerédi's regularity lemma to prove \emph{Roth's theorem on arithmetic progressions}, which states that every ``sufficiently dense'' set of natural numbers includes three elements of the form $k$, $k+d$, $k+2d$.
 
 We used a variety of source materials and discovered a good many significant infelicities in the definitions and proofs. These included confusion between $\subset$ and $\subseteq$ (which are often synonymous in combinatorics) and between a number of variants of the lemma statement. One minor claim was flatly incorrect. To make matters worse, the significance of these issues only became clear in the application of the regularity lemma to Roth's theorem. Much time was wasted, and yet the entire formalisation project~\cite{edmonds-formalising-regularity} took under six months.%
 \footnote{ An email from Angeliki proposing to prove Szemerédi's regularity lemma is dated 8 July 2021. The formalisation was done by 5 November; Roth, 28 December.}
By a remarkable coincidence, a group based in the mathematics department at Cambridge formalised a slightly different version of Szemerédi's regularity lemma, using Lean, around the same time~\cite{Dillies_Mehta_ITP_2022}.
 
 
\subsection{Additive combinatorics}

Let $A$ and $B$ be finite subsets of a given abelian group $(G,{+})$, and define their \emph{sumset} as
\[ A+B = \{a+b:a\in A, b\in B\}. \]
Write $nA$ for the $n$-fold iterated sumset $A+\cdots+A$.
\emph{Additive combinatorics} concerns itself with such matters as the relationship between the cardinality of $A+B$ and other properties of $A$ and~$B$. Angeliki proposed this field as the natural successor to the formalisation of Szemerédi's regularity lemma because it's fairly recent (many results are less than 50 years old) and significant (providing a route to Szemerédi's theorem, a much stronger version of the Roth result mentioned above).
 
 Here's an overview of the results formalised, all within the 7-month period from April to November 2022:
 \begin{itemize}
  \item The \emph{Plünnecke--Ruzsa inequality}: yields an upper bound on the \emph{difference} set $mB-nB$
  \item \emph{Khovanskii's theorem}: for any finite $A\subseteq G$, the cardinality of $nA$ grows like a polynomial for all sufficiently large $n$.
  \item The \emph{Balog–Szemerédi–Gowers theorem} is a deep result bearing on Szemerédi's theorem. The formalisation combines additive combinatorics with extremal graph theory and probability~\cite{koutsoukou-balog}.
  \item \emph{Kneser's theorem} and the \emph{Cauchy–Davenport theorem} yield lower bounds for the size of $A+B$.
\end{itemize}

These are highly significant results by leading mathematicians. They can all be found in Isabelle's \emph{Archive of Formal Proofs} (AFP).%
 \footnote{\url{https://www.isa-afp.org}} 

\subsection{Other formalisation projects}

The members chose a variety of large and small projects with a variety of specific objectives:
\begin{itemize}
  \item \emph{Combinatorial structures}. This is the PhD project of Chelsea Edmonds, who has used Isabelle's locale system to formalise dozens of varieties of block designs, hypergraphs, graphs and the relationships among them~\cite{edmonds-modular-first}. Results proved include Fisher's inequality~\cite{edmonds-formalising-fishers}.
  \item \emph{Number theory}. We have formalised several chapters of \emph{Modular Functions and Dirichlet Series in Number Theory}, a graduate textbook by Tom M. Apostol.
  \item \emph{Wetzel's problem} is a fascinating small example, due to Erdős, where the answer to a question concerning complex analysis depends on the truth or falsity of the continuum hypothesis. The formal proof illustrates analysis and axiomatic set theory smoothly combined into a single argument~\cite{paulson-wetzel-formalisation}.
  \item \emph{Turán's graph theorem} states a maximality property of Turán graphs. This was a Master's student project.
\end{itemize}

This is a partial list, especially as regards contributions from interns, students and other visitors.

\subsection{On legibility of formal proofs}

A proof is an argument, based on logical reasoning from agreed assumptions, that convinces mathematicians that a claim is true. How then do we understand a computer proof? To follow the analogy strictly, a computer proof convinces computers that a claim is true. But computers, even in this age of clever chatbots, are not sentient. We need to convince mathematicians.

Of the early efforts at the formalisation of mathematics, only Mizar aimed for legibility. Even pre-computer formal proofs such as \emph{Principia Mathematica} are unreadable. Isabelle's proof language (Isar) follows the Mizar tradition, as in the following example:

\begin{isabelle}
\isacommand{lemma}\ deriv\_sum\_int:\isanewline
\ \ "deriv\ (\isasymlambda x.\ \isasymSum i=0..n.\ real\_of\_int\ (c\ i)\ *\ x\isacharcircum i)\ x\ \isanewline
\ \ \ \ \ \ \ =\ (if\ n=0\ then\ 0\ else\ (\isasymSum i=0..n-1.\ of\_int((i+1)\ *\ c(Suc\ i))\ *\ x\isacharcircum i))"\isanewline
\ \ (\isakeyword{is}\ "deriv\ ?f\ x\ =\ (if\ n=0\ then\ 0\ else\ ?g)")\isanewline
\isacommand{proof}\ -\isanewline
\ \ \isacommand{have}\ "(?f\ has\_real\_derivative\ ?g)\ (at\ x)"\ \isakeyword{if}\ "n\ >\ 0"\isanewline
\ \ \isacommand{proof}\ -\isanewline
\ \ \ \ \isacommand{have}\ "(\isasymSum i\ =\ 0..n.\ i\ *\ x\ \isacharcircum \ (i\ -\ Suc\ 0)\ *\ (c\ i))\isanewline
\ \ \ \ \ \ \ \ =\ (\isasymSum i\ =\ 1..n.\ (real\ (i-1)\ +\ 1)\ *\ of\_int\ (c\ i)\ *\ x\ \isacharcircum \ (i-1))"\isanewline
\ \ \ \ \ \ \isacommand{using}\ that\ \isacommand{by}\ (auto\ simp:\ sum.atLeast\_Suc\_atMost\ intro!:\ sum.cong)\isanewline
\ \ \ \ \isacommand{also}\ \isacommand{have}\ "\isasymdots \ =\ sum\ ((\isasymlambda i.\ (real\ i\ +\ 1)\ *\ c\ (Suc\ i)\ *\ x\isacharcircum i)\ \isasymcirc \ (\isasymlambda n.\ n-1))\isanewline
\ \ \ \ \ \ \ \ \ \ \ \ \ \ \ \ \ \ \ \ \ \ \ \ \ \{1..Suc\ (n-1)\}"\isanewline
\ \ \ \ \ \ \isacommand{using}\ that\ \isacommand{by}\ simp\isanewline
\ \ \ \ \isacommand{also}\ \isacommand{have}\ "\isasymdots \ =\ ?g"\isanewline
\ \ \ \ \ \ \isacommand{by}\ (simp\ flip:\ sum.atLeast\_atMost\_pred\_shift\ [\isakeyword{where}\ m=0])\isanewline
\ \ \ \ \isacommand{finally}\ \isacommand{have}\ \isasymsection :\ "(\isasymSum a\ =\ 0..n.\ a\ *\ x\ \isacharcircum \ (a\ -\ Suc\ 0)\ *\ (c\ a))\ =\ ?g"\ \isacommand{.}\isanewline
\ \ \ \ \isacommand{show}\ ?thesis\isanewline
\ \ \ \ \ \ \isacommand{by}\ (rule\ derivative\_eq\_intros\ \isasymsection \ |\ simp)+\isanewline
\ \ \isacommand{qed}\isanewline
\ \ \isacommand{then}\ \isacommand{show}\ ?thesis\isanewline
\ \ \ \ \isacommand{by}\ (force\ intro:\ DERIV\_imp\_deriv)\isanewline
\isacommand{qed}
\end{isabelle}

Only a little training is required to make some sense of this. The lemma claims that the derivative of a certain summation equals a certain other summation. The proof refers of the variables \isa{?f} and \isa{?g}, which are defined by the pattern provided in the lemma statement: \isa{?f} denotes the original summation, and we prove that \isa{?g} is its derivative. Within that proof we can see summations being manipulated through changes of variable. Since we can see these details of the reasoning, we have reasons to believe that the proof is indeed correct: we do not simply have to trust the computer.

Not all Isabelle proofs can be written in a structured style. Page-long formulas often arise when trying to verify program code, and sometimes just from expanding mathematical definitions. Then we must use the traditional tactic style: long sequences of proof commands. However, most mathematical proofs that humans can write go into the structured style with ease. We have aimed for maximum legibility in all our work.

\section{Library Search and Machine Learning Experiments}

 The focus of this paper is achievements in the formalisation of mathematics, but the ALEXANDRIA proposal also called for investigating supporting technologies. The name of the project refers to the library of Alexandria, and Isabelle's AFP already has nearly 4 million lines of proof text and well over 700 separate entries. How can we take advantage of all this material when developing new proofs?

In May 2019, the team acquired a new postdoc: Yiannos Stathopoulos. He came with the perfect background to tackle these objectives. After much labour, he and Angeliki produced the SErAPIS search engine,%
\footnote{\url{https://behemoth.cl.cam.ac.uk/search/}}
which searches both the pre-installed Isabelle libraries and the AFP, offering a great many search strategies based on anything from simple keywords to abstract mathematical concepts~\cite{stathopoulos-developing}. It is not easy to determine the relevance or significance of a formal text to an abstract concept, but a variety of query types can be combined to explore the libraries. 

Also mentioned in the proposal was the aim of Intelligent User Support. I had imagined that common patterns of proofs could be identified in the existing libraries and offered up to users, but with no idea how. To generate structured proofs automatically would require the ability to generate intermediate mathematical assertions. Six years of dramatic advances in machine learning have transformed our prospects. Language models can generate plausible texts given a corpus of existing texts. And as the texts we want would be inserted into Isabelle proofs, we can immediately check their correctness. 

An enormous amount of work is underway, particularly by a student in our group, Albert Qiaochu Jiang, working alongside Wenda Li and others. It is now clear that language models can generate formal Isabelle proof skeletons~~\cite{li-isarstep} and can also be useful for identifying relevant lemmas~\cite{jiang-thor}. We can even envisage \emph{automatic formalisation}~\cite{jiang-draft-sketch,wu-autoformalization}: translating informal proofs into formal languages, by machine. 
 Autoformalisation is easier with a legible proof language like ours, because the formal proof can have the same overall structure as the given natural language proof; a project currently underway is to develop the Isabelle Parallel Corpus, pairing natural language and Isabelle texts.%
 \footnote{https://behemoth.cl.cam.ac.uk/ipc/}
 The next few years should see solid gains through machine learning.

\section{Evaluation}

At the start of this paper, I listed two scientific questions: what sort of mathematics, and what sort of proofs, can be formalised? And the answer so far is, everything we attempted, and we attempted a great variety of mathematical topics: number theory, combinatorics, analysis, set theory. The main difficulties have been errors and omissions in proofs. A vignette illustrates this point. Chelsea was formalising a probabilistic argument where the authors wrote ``these probabilities are clearly independent, and therefore the joint probability is obtained by multiplying them.'' The problem is that this multiplication law is the mathematical definition of independent probabilities, which the authors had somehow confused with the real-world concept of unconnected random events. Frequently we have found proofs that are almost right: they need a bit of adjustment, but getting everything to fit takes effort.

Effort remains the main obstacle to the use of verification tools by mathematicians. Obvious claims are often tiresome to prove, which is both discouraging and a waste of an expert's time. But we might already advocate an approach of formalising the definitions and the proofs, stating the obvious claims without proofs (using the keyword \textbf{sorry}). Even for this idea to be feasible, much more library material is needed,  covering at least all the definitions a mathematician might expect to have available.

Another key scientific question is the role of dependent types. People in the type theory world seem to share the conviction that dependent types are necessary to formalise nontrivial mathematics. But in reality it seems to be Lean users who repeatedly fall foul of  \emph{intensional equality}: that $i=j$ does not guarantee that $T(i)$ is the same type as $T(j)$. 
Falling foul of this can be fatal: the first definition of schemes had to be discarded for this reason.
Intensional equality is adopted by almost all dependent type theories, including Coq and Agda: without it, type checking becomes undecidable.
But with it, type dependence does not respect equality.

 The main limitation of simple type theory is that axiomatic type classes are less powerful than they otherwise would be. Isabelle/HOL has type classes for groups, rings, topological spaces among much else, but they are not useful for defining the theories of groups, rings or topological spaces. Rather they allow us, for example, to define the quaternions, prove a dozen or so laws and immediately inherit entire libraries of algebraic and topological properties. Abstract groups, rings, etc., need to be declared with an explicit carrier set (logically, the same thing as a predicate) rather than using the corresponding type class. It's a small price to pay for a working equality relation.
 
Having said this, one must acknowledge the enormous progress made by the Lean community over roughly the same period, 2017--now. Lean users, inspired by Buzzard, have taken on hugely ambitious tasks. The most striking is probably the Liquid Tensor Experiment~\cite{castelvecchi-mathematicians-welcome}: brand-new mathematics, by a Fields medallist (Peter Scholze) who was concerned about its correctness, formalised over about a year and a half. This one accomplishment, more than anything else, demonstrates that formalisation can already offer real value to professional mathematicians.

We have from time to time looked at type issues directly. De Vilhena~\cite{vilhena-algebraically-closed} describes an interesting technique for defining the $n$-ary direct product of a finite list of groups, iterating the binary direct product; his trick to avoid type issues involves creating an isomorphism to a suitable type. However, here one could avoid type issues (and handle the infinite case) by defining the direct product of a family in its own right as opposed to piggybacking off of the binary product. Anthony Bordg has done a lot of work on the right way to express mathematics without dependent types~\cite{bordg-encoding-dependently-typed,bordg-simple}. Ongoing work, still unpublished, is exploring the potential of the \textit{types-to-sets framework}~\cite{kuncar-from} to allow a smooth transition between type-based and carrier-set based formalisations.

 One can also compare formalisms in terms of their logical strength. Higher-order logic is somewhat weaker than Zermelo set theory, which is much weaker than ZFC, which in turn is much weaker than Tarski-Grothendieck set theory:
 \[ \mathrm{HOL} < \mathrm{Z} \ll \mathrm{ZF} \ll \mathrm{TG} \]
 The Calculus of Inductive Constructions, which is the formalism of Lean and Coq, is roughly equivalent to TG\@. The advantage of a weaker formalism is better automation. The power of ZF set theory, when it is required, can be obtained simply by loading the corresponding library from the AFP~\cite{paulson-wetzel-formalisation}. It's highly likely that a similar library could be created for Tarski-Grothendieck. And yet, remarkably, everything we have tried to formalise, unless it refers explicitly to ZF, sits comfortably within HOL alone. Since HOL is essentially the formalism of \emph{Principia Mathematica}~\cite{principia}, we can conclude that Whitehead and Russell were right all along.
 
The AFP entries contributed by the project authors are too many to list, but they can be consulted via the on-line author indices:
\begin{itemize}
  \item Anthony Bordg\\  \url{https://www.isa-afp.org/authors/bordg/}
  \item Chelsea Edmonds\\ \url{https://www.isa-afp.org/authors/edmonds/}
  \item Angeliki Koutsoukou-Argyraki\\  \url{https://www.isa-afp.org/authors/argyraki/}
  \item Wenda Li\\ \url{https://www.isa-afp.org/authors/li/}
  \item Lawrence C. Paulson\\ \url{https://www.isa-afp.org/authors/paulson/}
\end{itemize}

\section{Conclusions}

We set out to tackle serious mathematics with a combination of hope and trepidation. We were able to formalise everything we set out to formalise and were never forced to discard a development part way through. As Angeliki has pointed out, ``we have formalised results by two Fields medalists (Roth and Gowers), an Abel prize winner (Szemerédi) and of course Erdős too!''

 We've also seen impressive advances in search and language models to assist users in proof development. Although the effort required to formalise mathematical articles remains high, we can confidently predict that formalisation will be playing a significant role in mathematical research in the next few years.

\subsubsection{Acknowledgements} 
This work was supported by the ERC Advanced Grant ALEXANDRIA (Project GA 742178). Chelsea Edmonds, Angeliki Koutsoukou-Argyraki and Wenda Li provided numerous helpful comments and suggestions.

For the purpose of open access, the author has applied a Creative Commons Attribution (CC BY) licence to any Author Accepted Manuscript version arising from this submission.

\bibliographystyle{splncs04}
\bibliography{string,atp,general,isabelle,theory,crossref}
\end{document}